\PassOptionsToPackage{table}{xcolor}

\documentclass[a4paper]{article}

\usepackage{babel}
\usepackage[utf8]{inputenc}
\usepackage[margin=3.5cm]{geometry}
\usepackage{color}
\usepackage{graphicx}
\usepackage{xspace}
\usepackage{url}
\usepackage{amsmath,amssymb}
\usepackage{amsthm}
\usepackage{hyperref}
\usepackage[capitalize]{cleveref}
\usepackage{authblk}
\usepackage{todonotes}
\usepackage{paralist}
\usepackage{booktabs}
\usepackage{bbding}

\newtheorem{theorem}{Theorem}[section]

\newtheorem{problem}[theorem]{Problem}
\newtheorem{lemma}[theorem]{Lemma}

\newtheorem{proposition}[theorem]{Proposition}

\title{A tamed family of triangle-free graphs \\ with unbounded chromatic number}

\author{Édouard Bonnet}
\author{Romain Bourneuf}
\author{Julien Duron}
\author{Colin Geniet}
\author{Stéphan Thomassé}
\author{Nicolas Trotignon}
\affil{Univ.\ Lyon, ENS de Lyon, UCBL, CNRS, LIP, France}

\begin{document}

\maketitle

\begin{abstract}
We construct a hereditary class of triangle-free graphs with unbounded chromatic number,
in which every non-trivial graph either contains a pair of non-adjacent twins or has an edgeless vertex cutset of size at most two.
This answers in the negative a question of Chudnovsky, Penev, Scott, and Trotignon.
The class is the hereditary closure of a~family of (triangle-free) \emph{twincut graphs} $G_1, G_2, \ldots$ such that $G_k$ has chromatic number~$k$.
We also show that every twincut graph is edge-critical.
\end{abstract}

\section{Introduction}\label{sec:intro}

One of the main questions on the chromatic number $\chi(G)$ of a graph $G$ is how it compares to the clique number $\omega(G)$. Indeed, while $\omega(G) \leq \chi(G)$, early constructions by Blanche Descartes~\cite{BD54}, Zykov~\cite{ZY52}, and Mycielski~\cite{MY55} show that there are triangle-free graphs with arbitrarily large $\chi$. 
Such graphs have been an important source of inspiration in graph theory.
For instance, a distinctive early success of the probabilistic method was the construction by Erd\H{o}s~\cite{ER59} of graphs with large girth and large chromatic number. Another example is the proof by Lov\'asz~\cite{LO78} of the Kneser conjecture\footnote{Asserting that the Kneser graph $\mathcal K_{n,k}$, whose vertices are the $k$-subsets of $\{1,\ldots,n\}$ and whose edge relation is the disjointness of two sets, satifies $\chi(\mathcal K_{n,k})=n-2k+2$ for every $n \geqslant 2k$.}~\cite{Kneser55}, a cornerstone of the introduction of topological methods to combinatorial problems.

There is also an interesting interplay of these graphs with discrete
geometry in the plane.  For instance, triangle-free segment
intersection graphs were shown to have unbounded chromatic
number~\cite{PK13}, disproving a question of Erd\H{o}s and
Gy\'arf\'as~\cite{Gyarfas85}.  The proof consists of astutely
representing Burling graphs (another class of triangle-free graphs of
unbounded chromatic number that are intersection graphs of boxes of
$\mathbb R^3$)~\cite{BU65} as intersection graphs of segments in the
plane.  Recently, Davies~\cite{DA22} showed that the odd distance
graph on $\mathbb Z^2$ (with an edge between every pair of points at
Euclidean distance an odd integer) has infinite chromatic number, and
happens to be triangle-free, thereby providing another such class with
a geometric representation.

We build in this paper a new explicit sequence of triangle-free graphs
$G_k$, which we call \emph{twincut graphs}, satisfying $\chi(G_k)=k$
with the following striking property: all their induced subgraphs have
non-adjacent twins (two vertices with the same neighborhood), or an
edgeless vertex cutset of size at most two.  The details are given in
Section~\ref{sec:cons}.  This is very surprising since both situations
are, when considered individually, particularly favourable to keeping
the chromatic number low.  One the one hand, creating twins does not
change the chromatic number.  On the other hand, Alon, Kleitman, Saks,
Seymour and Thomassen \cite{AKSST87} proved that the closure of any
basic class under gluing along bounded subsets of vertices preserves
bounded chromatic number.  This was later refined by Penev, Thomassé
and Trotignon \cite{PTT16} who showed that such closure admits extreme
decompositions: a small vertex cutset isolates a basic subgraph of the
final graph, hence allowing a coloring with few colors.

A natural question is to consider two different types of closure, each
behaving well with respect to the chromatic number, and try to combine
them.  Along those lines, Chudnovsky, Penev, Scott, and Trotignon
\cite{CPST13} asked whether the closure of a $\chi$-bounded class
under substitutions and bounded cutsets could remain $\chi$-bounded,
where a \emph{$\chi$-bounded class} is a hereditary class of graphs
such that there exists a function $f$ satisfying
$\chi(G)\leq f(\omega(G))$ for all graphs of the class.

It may seem at first that this is just a matter of finding the right
induction hypothesis, but twincut graphs show that the answer is
negative in the seemingly easiest case: the closure $\cal C$ of
$\{K_1,K_2\}$ (the 1-vertex graph, and the edge) under the two
operations of vertex replication (i.e., creating a non-adjacent twin)
and gluing two graphs on up to two non-adjacent vertices.  To our
surprise, the class $\cal C$ turned out to contain all twincut graphs.

A~salient feature of constructions of large chromatic number is their
\emph{criticality}.  For instance, Kneser graphs are not
vertex-critical (their chromatic number need not drop when a vertex is
removed), hence Schrijver \cite{SC78} proposed a canonical way to
pinpoint a critical induced subgraph with the same chromatic number.
There are very few critical constructions\footnote{In a strict
  explicit and deterministic sense, since one can always greedily
  remove edges.} and to our knowledge, only the Mycielski sequence and
its generalized variants achieve edge-criticality.  In terms of
structural complexity, Mycielski graphs are universal (they contain
all triangle-free graphs as induced subgraphs), and their generalized
counterparts have unbounded Vapnik-Chervonenkis dimension (they
contain all bipartite graphs as induced subgraphs).  Surprisingly,
twincut graphs achieve edge-criticality while keeping low VC-dimension
(for example, they do not induce the cube).

In a forthcoming paper, we compute several
width-parameter values (tree-width, rank-width and twin-width) of
twincut graphs. This confirms their very low structural complexity.
We also provide a full structural description of
the class formed by the induced subgraphs of the graphs $G_k$ together
with a polynomial time recognition algorithm and evidence that twincut
graphs are related to previous constructions (namely that every
twincut graph is an induced subgraph of some Zykov graph and a
(non-induced) subgraph of some Burling graph).

\section{The twincut graphs}
\label{sec:cons}

A \emph{structured tree} is a rooted tree~$T$ and a function~$g$ defined on the internal nodes~$v$ of $T$ (i.e. non leaves) such that~$g(v)$ is a graph whose vertices are the children of~$v$ in~$T$. A \emph{branch} in $T$ is a path from the root to one of the leaves of $T$. 
The \emph{realization} $R(T,g)$ of $(T,g)$ is the graph defined on vertex set $V(T)\cup B$, where $B$ is the set of branches of $T$. The edges of $R(T,g)$
first consist of all $uv$ where $u,v$ are children of $z$ and $uv$ is an edge of $g(z)$.
At this point, the graph $R(T,g)$ is simply the disjoint union of all $g(z)$ and some isolated vertices ($B$ and the root).
Next, we connect each \emph{branch vertex} $b\in B$ to all the vertices of $T$ in the branch $b$. Observe that the edges of $T$ are not edges of $R(T,g)$.

Note that when $T$ has only one (root) vertex, it is also a leaf. In particular $g$ is empty ($T$ has no internal node) and therefore $R(T, g)$ is obtained from $T$ by adding a single vertex which is adjacent to the root. Hence, $R(T, g)$ is $K_2$.

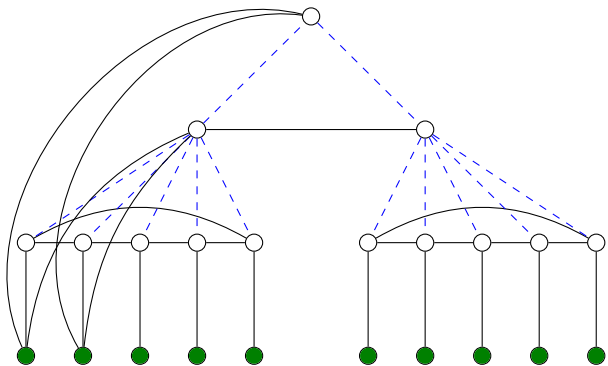
\begin{figure}[h!]
\centering
\begin{tikzpicture}[vertex/.style={draw,circle,inner sep=-0.08cm},leaf/.style={fill=green!50!black,circle,inner sep=-0.075cm},edget/.style={thin,blue,dashed}]
\def\s{1.5}

\node[vertex] (t) at (0,3 * \s) {} ;

\foreach \z/\l in {-1/r,1/l}{
\foreach \i/\j/\k in {-1/2/t,-2.5/1/a1,-2/1/a2,-1.5/1/a3,-1/1/a4,-0.5/1/a5,
-2.5/0/b1,-2/0/b2,-1.5/0/b3,-1/0/b4,-0.5/0/b5}{
\node[vertex] (\k\l) at (\i * \s * \z,\j * \s) {};
}
\foreach \i in {1,...,5}{
\node[leaf] at (b\i\l) {};
}

\foreach \i [count=\ip from 2] in {1,...,4}{
\draw (a\i\l) -- (a\ip\l) ;
\draw (a\i\l) -- (b\i\l) ;
}
\draw (a1\l) to [bend left=30 * \z] (a5\l) ;
\draw (a5\l) -- (b5\l) ;
\foreach \i in {1,...,5}{
\draw[edget] (t\l) -- (a\i\l) ;
}
}

\draw (tl) -- (tr) ;

\draw[edget] (tl) -- (t) -- (tr) ;

\draw (b1l) to [bend left=30] (tl) ;
\draw (b1l) to [bend left=64] (t) ;
\draw (b2l) to [bend left=20] (tl) ;
\draw (b2l) to [bend left=64] (t) ;

\end{tikzpicture}
\caption{The 4-chromatic triangle-free graph $G_4$. The tree $T_4$ is represented with dashed blue edges (which are \emph{not} actual edges of $G_4$). Every green vertex is adjacent to all vertices in a branch of $T_4$. We explicitly represented these edges for the two leftmost green vertices.}
\label{fig:G4}
\end{figure}

We present now an inductive construction of a family of triangle-free graphs $(G_i)_{i \in \mathbb N^+}$, called \emph{twincut graphs}, with unbounded chromatic number.
First,~$G_1$ is defined as the graph on one vertex.
Assuming that $G_1,\dots,G_{k-1}$ have been built, the graph~$G_k$ is defined as the realization of the following structured tree~$(T_k, g_k)$:
the tree~$T_k$ has~$k-1$ levels (the root being at level 1), and for each node~$v$ at level~$i<k-1$, we give~$|V(G_{i+1})|$ children to~$v$ and set $g_k(v) = G_{i+1}$. For instance~$T_2$ consists only of its root, and its realization~$G_2$ is $K_2$ as explained above. Then,~$T_3$ has a root~$r$ with two children~$c,c'$ which are linked in~$g_3(r) = G_2$.
The realization adds a vertex~$x$ connected to~$r,c$, and a vertex $y$ connected to~$r,c'$, thus creating a 5-cycle $rxcc'y$, hence $G_3 = C_5$.
The graph $G_4$ has $1+2+10+10=23$ vertices, see \cref{fig:G4}.

\begin{proposition}\label{prop:triangle-free}
For every integer $k \geq 1$, $G_k$ is triangle-free.
\end{proposition}
\begin{proof}
This can be seen by induction on $k$.
$G_1$ is triangle-free since it has a single vertex.
$G_{k+1}$ is obtained from the disjoint union of copies of $G_1, G_2, \ldots, G_k$, which by the induction hypothesis is triangle-free, by adding vertices adjacent to an independent set.
Indeed each new vertex $b$ in $G_{k+1}$ is adjacent to at most one vertex in each copy of the graphs $G_1, G_2, \ldots, G_k$, hence cannot create a triangle.
Thus $G_{k+1}$ is itself triangle-free.
\end{proof}

Twincut graphs have unbounded chromatic number, with a~similar argument to the one used for Zykov graphs, and the additional twist of finding a~rainbow independent set along a branch of the structured tree.
\begin{proposition}\label{prop:high-chi}
For every integer $k \geq 1$, we have $\chi(G_k) = k$.
\end{proposition}
\begin{proof} 
The proof is again by induction on $k$. 
The case $k=1$ holds since $G_1$ is a 1-vertex graph.
Now, let $k \geq 1$ and suppose $\chi(G_\ell) = \ell$ for $\ell \leq k$.
Fix $c$ a proper coloring of $G_{k+1}$.
In the underlying structured tree~$T_{k+1}$, we will pick a branch which uses~$k$ distinct colors.
Assume by induction that~$v_1,\dots,v_\ell$ is a path in~$T_{k+1}$ starting from the root~$v_1$ such that the colors~$c(v_i)$ are all distinct.
By construction of~$G_{k+1}$, the children of~$v_\ell$ induce a copy of~$G_{\ell+1}$,
which is $(\ell+1)$-chromatic.
Thus, there is a child~$v_{\ell+1}$ whose color is distinct from~$c(v_1),\dots,c(v_\ell)$,
with which we extend the path.
Once this process reaches a~leaf of~$T_{k+1}$, we obtain a branch~$b$ whose vertices use~$k$ distinct colors,
hence the vertex~$b$, which is connected exactly to this branch, needs one additional color.
Thus, $c$ uses at least $k+1$ colors, so $\chi(G_{k+1}) \geq k+1$.

Conversely, if we color in~$G_{k+1}$ all branch vertices by~$k+1$ and remove them from $G_{k+1}$,
we are left with the disjoint union of all graphs~$g(v)$, i.e., copies of~$G_1,\dots,G_k$ which are $k$-colorable by induction.
This yields a $(k+1)$-coloring of~$G_{k+1}$.
\end{proof}

In \cite{CPST13}, the authors show that the closure of a $\chi$-bounded class  under substitution is~\mbox{$\chi$-bounded}, and that substitutions further preserve polynomial $\chi$-boundedness. This is also true when the closure consists of gluing pairs of graphs  along bounded size subsets. Trying to merge these two operations, they posed the following problem, also mentioned in~\cite{SS20}:

\begin{problem}\label{prob:main}
Is the closure of a $\chi$-bounded class under substitution and gluing along a~bounded
number of vertices also $\chi$-bounded?
\end{problem}

Twincut graphs give a strong negative answer to \cref{prob:main}.
Let~$\mathcal{C}$ be the closure of the graphs of size at most two under the following two operations: substituting a vertex by a stable set of size two, and gluing two graphs of~$\mathcal{C}$ along a stable set of size at most two. This definition is a very special case of the closure considered in \cref{prob:main}. Observe that the class~$\mathcal{C}$ is closed under taking induced subgraphs. Note also that the  graphs in~$\mathcal{C}$ are triangle-free.
Thus, to negatively answer \cref{prob:main} it suffices to prove the following:

\begin{proposition}
    \label{prop:Tr-twincut}
    The graphs~$G_k$ are in~$\mathcal{C}$.
\end{proposition}
We more generally show that~$\mathcal{C}$ is closed under the realization of structured trees,
which immediately implies \cref{prop:Tr-twincut}.
\begin{lemma}
    \label{lem:Tr-realization}
    Let~$(T,g)$ be a structured tree such that every~$g(v)$ is in~$\mathcal{C}$. Then $R(T,g) \in \mathcal{C}$.
\end{lemma}
\begin{proof}
    For a node~$v$ of~$T$, let~$T(v)$ be the subtree rooted at~$v$, i.e., the subtree consisting of all descendants of~$v$. Equipped with the restriction of~$g$, $T(v)$ is a structured tree.
    We prove by induction on~$T$, starting from the leaves, that for all nodes~$v$, the realization~$R(T(v), g)$ is in~$\mathcal{C}$.
    For the sake of brevity, let us denote this realization of a subtree by~$R_v = R(T(v),g)$.

    If~$v$ is a leaf, then~$R_v$ is simply an edge, which is in~$\mathcal{C}$.
    Let now~$v$ be an internal node with children~$u_1,\dots,u_\ell$, and assume that each~$R_{u_i}$ is in~$\mathcal{C}$.
    Recall also that~$g(v)$ is assumed to be in~$\mathcal{C}$.
    We construct~$R_v$ as follows. First, in each~$R_{u_i}$, create a copy~$u'_i$ of~$u_i$ by substituting~$u_i$ with a stable set of size~two, and call~$R'_{u_i}$ the resulting graph.
    Next, take~$g(v)$ and add to it an isolated vertex standing for~$v$.
    We then glue each~$R'_{u_i}$ successively with this graph,
    by identifying~$u'_i$ with~$v$, and identifying the occurrences of~$u_i$ in~$R'_{u_i}$ and in~$g(v)$.
    This corresponds to gluing along a stable set of size~two.
    Hence, we constructed~$R_v$ starting from~$g(v),R_{u_1},\dots,R_{u_\ell}$, by substituting with and gluing on stable sets of size at most~two, thereby proving that~$R_v \in \mathcal{C}$.
\end{proof}

\section{Criticality of twincut graphs}
\label{sec:min}

Recall that a graph $G$ is \emph{critical} (or \emph{edge-critical}) if every strict subgraph $H$ of $G$ satisfies $\chi(H)<\chi(G)$. In other words, deleting an edge from $G$ decreases its chromatic number.
\begin{proposition} \label{prop:critical}The graphs $G_k$ are critical.
\end{proposition}

\begin{lemma} For every $k \geq 1$, for every vertex $v$ of $G_k$, there exists a proper $k$-coloring of $G_k$ in which $v$ is the only vertex with color $k$. Furthermore, if $v \in B$ then for every $i$, the vertex of $N(v)$ at level $i$ in $T_k$ has color $i$.
\end{lemma}

\begin{proof}
    The proof is by induction on $k$. The property holds for $k = 1$. Let $k \geq 1$ and assume that the property holds for every $\ell \leq k$. Let $v$ be any vertex of $G_{k+1}$. 
    \begin{itemize}
        \item If $v$ is a branch vertex $b$ where $b = v_1, v_2, \ldots, v_k$, then color the root by $c(v_1) = 1$, and for every $j > 1$, fix a proper $j$-coloring of $g_{k+1}(v_{j-1})$ in which $v_j$ is the only vertex of color $j$. For every other vertex $w$ in $T_{k+1}$ at some level $\ell < k$, fix an arbitrary $(\ell+1)$-coloring of $g_{k+1}(w)$. At this point, every vertex of $G_{k+1}$ has a color in $\{1, \ldots, k\}$, except for the branch vertices. Let $b' \neq b$ be a branch of $T_k$. Write $b' = w_1, \ldots, w_k$. Take $\ell$ minimum such that $w_\ell \neq v_\ell$ (note that $\ell > 1$). Then, $w_\ell$ is a vertex of $g_{k+1}(v_{\ell-1})$ so $w_\ell$ does not have color $\ell$ by definition of the coloring of $g_{k+1}(v_{\ell-1})$. Hence, not all colors in $\{1, \ldots, k\}$ appear in the branch $b'$. Thus, we can color the branch vertex $b'$ with some color from $\{1, \ldots, k\}$. Finally, set $c(v) = k+1$.
        \item If $v$ is a vertex of $T_{k+1}$, pick an arbitrary branch $b$ of $T_k$ containing $v$. Like above, color $G_{k+1}$ so that the branch vertex $b$ is the only vertex of color $k+1$, and all its neighbors have different colors. Finally, swap the colors of~$v$ and~$b$.
    \end{itemize}
\end{proof}

We can now prove Proposition \ref{prop:critical}.

\begin{proof}
Let $uv$ be an edge in $G_k$. Let us show that $G_k\setminus uv$ is $(k-1)$-colorable.
\begin{itemize}
    \item If $u$ is a branch vertex $b$ of $G_k$, consider a proper $k$-coloring of $G_{k}$ in which $u$ is the only vertex of color $k$. In $G_k\setminus uv$, $u$ has degree $k-2$ so we can recolor it with some color from $\{1, \ldots, k-1\}$. The same holds if $v$ is a branch vertex $b$.
    \item If both $u,v$ belong to some graph $g(w)$ at level $i$ in $T_k$. Fix a branch $b$ containing $w$ and consider a proper $k$-coloring of $G_{k}$ in which the branch vertex $b$ is the only vertex of color $k$ and such that for every $j$, the vertex of the branch $b$ at level $j$ has color $j$. Since $g(w)$ is $i$-critical, we can recolor $g(w)$ using colors $\{1,\dots ,i-1\}$. We can then recolor the branch vertex $b$ with color $i$. 
\end{itemize}
\end{proof}

To our knowledge, the only explicit construction of critical high chromatic triangle-free graphs is the sequence of (generalized) Mycielski graphs. This class has high complexity since Mycielski graphs contains \emph{all} triangle-free graphs as induced subgraphs. To pinpoint a relevant complexity measure, we can note that the Vapnik-Chervonenkis dimension of the class of Mycielski graphs is unbounded. Stated in a less formal (albeit equivalent) way: all bipartite graphs appear as induced subgraphs of Mycielski graphs. This is also the case in their generalized version with large odd girth.

Another classical construction, the Zykov graphs, also have unbounded VC-dimension since the $k^{th}$ iteration already contains all bipartite graphs of size $(k-2,k-2)$. Twincut graphs form a subclass of Zykov graphs with bounded VC-dimension. Indeed, the cube is not an induced subgraph of a twincut graph: it has no twins, and no vertex cutset of size at most 2.


\end{document}